\documentclass{amsart}
\usepackage{amscd,amssymb}
\theoremstyle{plain}
\newtheorem{prop}[subsection]{Proposition}
\newtheorem{thm}[subsection]{Theorem}
\newtheorem{lem}[subsection]{Lemma}

\theoremstyle{remark}

\theoremstyle{definition}

\numberwithin{equation}{section}

\newcommand{\bF}{{\mathbb F}}

\newcommand{\Z}{{\mathbb Z}}
\newcommand{\C}{{\mathbb C}}

\begin{document}

\title[Orbit configuration spaces associated to discrete subgroups
of $PSL(2,\mathbb R)$]
{Orbit configuration spaces associated to discrete subgroups
of $PSL(2,\mathbb R)$}
\author{F. R.~Cohen*}
\address{Department of Mathematics \\
University of Rochester\\ Rochester, NY 14627 U.S.A.}
\email{cohf@math.rochester.edu}
\thanks{*Partially supported by the NSF grant number 9704410}

\author{Toshitake KOHNO}
\address{Department of Mathematical Sciences\\ University of
Tokyo\\ Tokyo 153-8914 Japan}
\email{kohno@ms.u-tokyo.ac.jp}

\author{M. A.~Xicot\'encatl**}
\address{Depto. de Matem\'aticas \\
Cinvestav \\ Apdo. Postal 14-740 \\
M\'exico, D.F. 07300.}
\email{xico@math.cinvestav.mx}
\thanks{**Partially supported by CONACyT grant number 37296-E}


\begin{abstract}

The purpose of this article is to analyze several Lie algebras
associated to ``orbit configuration spaces'' obtained from a group
$G$ acting freely, and properly discontinuously on the upper $1/2$-plane
$\mathbb H^2$. The Lie algebra obtained from the descending central
series for the associated fundamental group is shown to be
isomorphic, up to a regrading, to

\begin{enumerate}
    \item the Lie algebra obtained from
the higher homotopy groups of ``higher dimensional arrangements''
modulo torsion, as well as
    \item the Lie obtained from horizontal chord diagrams for
    surfaces.
\end{enumerate} The resulting Lie algebras are similar to those
studied in \cite{K1,K2,K3,C,CX,FH,CS}.
The structure of a related graded Poisson algebra defined
below and obtained from an analogue of the infinitesimal braid relations
parametrized by
$G$ is also addressed.
\end{abstract}

\maketitle

\section{Introduction}
The purpose of this article is to consider
certain Lie algebras which are described next.
\begin{enumerate}
\item One Lie algebra is obtained from the descending central series for
the fundamental group of orbit configuration spaces
associated to discrete subgroups of $PSL(2,\mathbb R)$
acting freely, and properly discontinuously on the upper $1/2$-plane
$\mathbb H^2$
by fractional linear transformations.
\item  A second Lie algebra is obtained from the classical higher
homotopy groups modulo torsion for loop spaces of orbit configuration
spaces of points in $\mathbb H^2 \times \mathbb C^q$.
\item A third Lie algebra is obtained from
horizontal chord diagrams for surfaces.
\end{enumerate}The main results here are that these Lie algebras,
apart from a ``trivial" degree shift, are
isomorphic. In addition, there are natural Poisson algebras obtained
from the homology of the iterated loop spaces of these orbit
configuration spaces, a structure given below.

Orbit configuration spaces were studied by the third author in \cite{X}, and
are defined next. Given a manifold $M$ on which the discrete group $G$ acts
properly
discontinuously, let ${\rm Conf}^{G}(M, n)$ denote the orbit
configuration space  $${\rm Conf}^{G}(M,n) =
\{(m_1, \dots, m_n) \in M^n \mid G \cdot m_i \cap G \cdot m_j =
\varnothing \text{\; if \; } i\neq j\}.$$ It will be assumed
throughout this article that $M \to\ M/G$ is the projection
map for a covering space.

If $G$ is a discrete group, then by \cite {X}, there are
fibrations $${\rm {\rm Conf}}^{G}(M,n) \to\ {\rm {\rm
Conf}}^{G}(M,i)$$ with fibre over the point
$(p_1,p_2,\cdots,p_i)$ in ${\rm Conf}^{G}(M,i)$ given by
$${\rm Conf}^{G}(M-Q_i^{G},n-i)$$ where $G \cdot p$
denotes the $G$-orbit of $p$, and $$Q_i^{G}= \amalg_{1
\leq j \leq i} G \cdot p_j.$$  The groups, and manifolds
addressed in this article are given as follows.
\begin{enumerate}
\item The space $M$ is either the upper $1/2$-plane $\mathbb H^2 =
SL(2, \mathbb R )/SO(2)$ or a product $\mathbb H^2 \times \mathbb C^q$.
\item The group $G$ is a discrete subgroup of $PSL(2, \mathbb R )$ where it
is assumed that $G$ acts freely on $\mathbb H^2 $ by
fractional linear transformations, trivially on $\mathbb C^q$, and
diagonally on the product $\mathbb H^2 \times \mathbb C^q$. In particular,
the natural projection $$\mathbb H^2 \to\ \mathbb H^2 / G$$ is
the projection for a covering space.
\end{enumerate}

Next consider the Lie algebra obtained from the descending central
series for a discrete group $G$. For each strictly positive integer $q$,
there is a canonical (and trivially defined) graded Lie algebra
$E_0^*(G)_q$ attached to the one obtained from the descending
central series for $G$, and which is defined as follows.

\begin{enumerate}
\item Fix a strictly positive integer $q$.
\item Let $\Gamma^m(G)$ denote the $m$-th
stage of the descending central series for $G$.
\item $E_0^{2mq}(G)_q = \Gamma^m(G) /  \Gamma^{m+1}(G),$
\item $E_0^{i}(G)_q = \{0\}$, if $i \not\equiv 0 \mod 2q$,  and
\item the Lie bracket is induced by that for the associated
graded for the $\Gamma^m(G)$.
\end{enumerate}

Restrict attention to path-connected topological spaces $X$ which
either have torsion free homology or where homology is taken with field
coefficients. This technical condition implies the strong form of the
K\"unneth theorem gives that the homology of $X \times X$ is isomorphic
to $H_*X \otimes H_*X$. In this case, the module of primitive elements
in the homology of $X$, $\operatorname{Prim}H_*(X)$, is the module
generated by elements
$\alpha$
which have trivial coproduct. The universal enveloping algebra of a graded
Lie algebra $L$ is denoted $U[L]$. Fix a base-point
${\mathbf x}=(x_1, \cdots, x_n)$ in ${\rm {\rm Conf}}(S_g,n)$. Let
${\mathbf y}$ denote a fixed choice of
base-point for ${\rm Conf}^{G}(\mathbb H^2,n)$, and $x_0$ a base-point for
$S_g$. The fundamental group of the surface $S_g$ is
denoted $\pi_1(S_g,x_0) = G$,
and the fundamental group of the orbit configuration space is denoted
$\pi_1({\rm Conf}^{G}(\mathbb H^2,n), {\mathbf y}) = P_n(S_g)^0$.

\begin{thm}\label{thm:Lie algebras}
Let $n$, and $q$ be fixed natural numbers.
There are isomorphisms of Lie algebras
$$E^*_0( P_n(S_g)^0 )_q \to\ \operatorname{Prim}H_*\Omega
({\rm Conf}^{G}(\mathbb H^2 \times \mathbb C^q,n)),$$ and
$$U[E^*_0(P_n(S_g)^0)_q] \to\ H_*\Omega
({\rm Conf}^{G}(\mathbb H^2 \times \mathbb C^q,n)).$$ Furthermore,
these Lie algebras are prescribed as follows:
\begin{enumerate}
            \item There are sub-Lie algebras
            of $E^*_0(P_n(S_g)^0)_q$
given by $L[i]$ the free Lie algebra generated by elements
$B_{i,j}^{\sigma}$ for fixed i
with $ 1 \leq j < i \leq n $ of degree $2q$ and where $\sigma$ runs over
the elements of the
group  ${G}$.
          \item There is an isomorphism of abelian groups given by the
natural additive extension $$\bigoplus_{ 2 \leq i \leq n} L[i] \to\
E^*_0(P_n(S_g)^0)_q.$$
          \item A complete set of relations are as follows.
              \begin{enumerate}
            \item If
            $\{i,j\} \cap \{s,t\} = \phi$, then $[B_{i,j}^{\sigma},
            B_{s,t}^{\tau}] = 0$.
            \item If $ 1 \leq j < s < i \leq k$, then
            $[B^{\tau}_{i,j},B^{\tau \sigma^{-1} }_{i,s}+ B^{\sigma}_{s,j}] =
            0$.
            \item If $ 1 \leq j < s < i \leq k $, then
            $[B^{\sigma}_{s,j},B^{\tau }_{i,j} + B^{\tau \sigma^{-1}}_{i,s}] =
            0$.
            \item The antisymmetry relation, and Jacobi identity for a graded
            Lie algebra.
            \end{enumerate}
\end{enumerate}
\end{thm}
\medskip

The symmetric group on $n$-letters $\Sigma_n$ acts naturally on
the configuration space, and thus on the cohomology ring. On the
other hand, the symbol $B^{\gamma }_{j,i}$ has not been defined
above in the cases for which $j<i$. This element is defined to
be $\tau(i,j)(B^{\gamma}_{i,j})$ where $\tau(i,j)$ is the
permutation which switches $i$, and $j$. The element
$\tau(i,j)(B^{\gamma}_{i,j})$ will
be shown to be equal to $B^{{\gamma}^{-1} }_{i,j}$ in
Lemma \ref{lem:Lie antisymmetry} of section $3$ below. In addition
the full action of the symmetric group is also specified
in Lemma \ref{lem:Lie antisymmetry}. The
following additional properties are satisfied
for ${\rm Conf}^{G}(\mathbb H^2,n)$.

\begin{thm}\label{thm:the case H}
Let $n$ be fixed natural number.
The following properties are satisfied.
\begin{enumerate}
\item The orbit configuration space ${\rm Conf}^{G}(\mathbb H^2,n)$ is a
$K(P_n(S_g)^0,1)$.
The wreath product of the
symmetric group $\Sigma_n$ with $G$, $\Sigma_n \wr G$,
acts properly discontinuously on ${\rm Conf}^{G}(\mathbb H^2,n)$.
The orbit space $${\rm Conf}^{G}(\mathbb H^2,n) / \Sigma_n \wr G$$ is
homeomorphic to ${\rm Conf}( \mathbb H^2 / G,n)/\Sigma_n $, a $K( \Pi,1)$
where $\Pi$ is the $n$-stranded braid group for the surface $\mathbb H^2/G$.

\item The fibration ${\rm Conf}^{G}(\mathbb H^2,n) \to\ {\rm
Conf}^{G}(\mathbb H^2,n-1)$
has (i) trivial local coefficients in homology, and (ii) a cross-section.

\item The Lie algebra $E^*_0( \pi_1({\rm Conf}^{G}(\mathbb H^2,n))$
is isomorphic to that given in Theorem \ref{thm:Lie algebras}.

\item The integral homology of ${\rm Conf}^{G}(\mathbb H^2,n)$ is given
additively by
$$ H_*{\rm Conf}^{G}(\mathbb H^2,n) \cong  H_* (C_1) \otimes H_* (C_2) \otimes
\dots \otimes H_* (C_{n-1})$$ where $C_i$ is the infinite bouquet of
circles $\bigvee_{|Q_i^{G}|} S^1$ where $Q_i^{G}$ as defined
as above.
\end{enumerate}
\end{thm}
\medskip

Horizontal chord diagrams of a closed oriented
surface of genus $g$, $S_g$, are described in section $2$ of this
article. In case $ g  > 0 $, these horizontal diagrams give
analogous constructions as those in genus $0$ as given in
\cite{K2,K3,K4}. Let ${\mathcal A}_n(S_g)^0$ and ${\mathcal
A}_n(S_g)$ denote the algebras of horizontal chord
diagrams of $S_g$  as introduced, and studied in
\cite{GP} (see section 2 for the notation).

\begin{thm}\label{thm:chord diagrams}
The algebra $${\mathcal A}_n(S_g)^0$$
is  isomorphic to the universal enveloping algebra
of $E^*_0( \pi_1{\rm Conf}^{G}(\mathbb H^2,n),\mathbf y) =
E^*_0(P_n(S_g)^0)$
where $$S_g = \mathbb H^2/G.$$
\end{thm}
\medskip

Work of T. Kohno \cite{K1, K2,K3}, as well as subsequent work of
M. Falk and R. Randell \cite{FR} gives the structure of the Lie
algebra associated to the descending central series for the
$k$-stranded pure braid group, the fundamental group of ${\rm
Conf}(\mathbb R^2,n)$. Work in \cite{CX} gives the associated Lie
algebras for loop spaces of
${\rm Conf}^{L}(\mathbb C \times {\mathbb C}^q,n)$ where $L$
is the the standard lattice of integral points in $\mathbb C$ acting
by translations on $\mathbb C$. Work of J.~Gonz\'alez-Meneses,
and L.~Paris \cite{GP} gives the
structure of the Lie algebras associated to the descending central
series for the fundamental group of ${\rm Conf}^{G}(\mathbb H^2,n)$.

These Lie algebras arise in an analysis of Vassiliev
invariants for pure braids on a surface. The theorem above gives a
comparison between these Lie algebras, and their analogues for
certain choices of loop spaces. Analogous groups were shown to be
given by morphisms of coalgebras in \cite{CS}. In addition, recent
work of Papadima, and Suciu \cite{PS} give related structures when
the associated Lie algebras are of finite type. The ones which
occur here for a surface of genus greater than $0$ are not of
finite type, and it is as yet unclear whether there is an
extension of the work in \cite{PS} to these choices of Lie
algebras.

These Lie algebras also occur in a second convenient context.
Namely, by \cite{X} there is a principal $G^n$ bundle
$${\rm Conf}^{G}(\mathbb H^2,n) \to\ {\rm Conf} (\mathbb H^2/ G,n), $$ and thus
there is a short exact sequence of groups
$$1 \to\ \pi_1({\rm Conf}^{G}(\mathbb H^2,n),{\mathbf y})
\to\  \pi_1({\rm Conf}(\mathbb H^2/ G,n),{\mathbf x})
\to\ G^n\to\ 1.$$ Hence the kernel of the map
$$ \pi_1({\rm Conf}(\mathbb H^2/ G,n),{\mathbf x})
\to\ G^n $$ induced by the natural inclusion
$$ {\rm Conf}(\mathbb H^2/ G,n)
\to\ (\mathbb H^2/ G)^n $$ is isomorphic to the fundamental group
$\pi_1({\rm Conf}^{G}(\mathbb H^2,n),{\mathbf y}) = P_n(S_g)^0$. The
Lie algebra obtained from the descending
central series for $\pi_1(\rm Conf(\mathbb H^2/ G,n),{\mathbf x})$ is not
analyzed here. The structure of the analogous Lie algebra
associated to $\pi_1({\rm Conf}^{G}(\mathbb H^2,n),{\mathbf y})$
is handled more easily than that for
$\pi_1({\rm Conf}(\mathbb H^2/ G,n),{\mathbf
x})$ and is given above, as well as earlier in \cite{GP}.

There is a natural structure of graded
Poisson algebra structure for the homology
of an iterated loop space. This Poisson algebra given by
$$H_*\Omega^k ({\rm Conf}^{G}(\mathbb H^2 \times \mathbb
C^q,n))$$ which is described in more detail in section $6$ here. Loosely
speaking, the resulting structure is the free Poisson algebra subject to
the ``Poisson analogues" of the relations in Theorem \ref{thm:Lie
algebras}. Precise definitions in the next theorem are given in
section $4$.

\begin{thm} \label{thm:poisson}
Assume that $k$ is least $1$.
\begin{enumerate}
\item If $k > 1$, the homology of
$\Omega^k({\rm Conf}^{G}(\mathbb H^2 \times \mathbb C^q,n))$, with any
field coefficients, is a graded Poisson algebra with Poisson bracket
given by the Browder operation $\lambda_{k-1}[-,-]$
for the homology of a $k$-fold loop space.

\smallskip

\item If $1< k <2q+1$, the homology of
$\Omega^k {\rm Conf}^{G}(\mathbb H^2 \times \mathbb C^q,n)$
with coefficients in a field $\bF$ of characteristic zero, is
the free Poisson algebra generated by elements
$$B_{i,j}^{\sigma}$$ of degree $2q+1 -k$ for $1 \leq j < i \leq n$, and
$\sigma$ in
$G$ modulo the ``infinitesimal Poisson surface
braid relations" given as follows:

      \begin{enumerate}
            \item If
            $\{i,j\} \cap \{s,t\} = \phi$, then
            $\lambda_{k-1}[B_{i,j}^{\sigma}, B_{s,t}^{\tau}] = 0$.
            \item If $ 1 \leq j < s < i \leq n$, then
            $\lambda_{k-1}[B^{\tau}_{i,j},B^{\tau \sigma^{-1} }_{i,s}+
            B^{\sigma}_{s,j}] = 0$.
            \item If $ 1 \leq j < s < i \leq n $, then
            $\lambda_{k-1}[B^{\sigma}_{s,j},B^{\tau }_{i,j} + B^{\tau
\sigma^{-1}}_{i,s}] =0$.
            \item The antisymmetry relation, and Jacobi identity for a graded
Poisson algebra.
      \end{enumerate}

\smallskip

      \item There is a map $$E^2: {\rm Conf}^{G}(\mathbb H^2 \times \mathbb
C^q,n) \to\ \Omega^2{\rm Conf}^{G}(\mathbb H^2 \times \mathbb
C^{q+1},n) $$ which induces a homology isomorphism in degree $2q
+1$. The associated loop map $\Omega(E^2): \Omega {\rm
Conf}^{G}(\mathbb H^2 \times \mathbb C^q,n) \to\ \Omega^3{\rm
Conf}^{G}(\mathbb H^2 \times \mathbb C^{q+1},n)$ induces an
isomorphism on $H_{2q}(-;\Z)$. Furthermore, the image of the map
$\Omega(E^2)$ in homology is the subalgebra generated by the
classes of degree $2q$.

\end{enumerate}
\end{thm}
\medskip

\section{Horizontal chord diagrams for surfaces}

Let $S_g$ be a closed oriented surface of genus $g$.
The main subject of this section is a relation between
Vassiliev invariants for
the braid groups of $S_g$ and
horizontal chord diagrams on $S_g$.
We refer the reader to \cite{BN}
for Vassiliev invariants of knots in $S^3$
and chord diagrams.
Let ${\rm Conf}(S_g,n)$ denote the configuration space of
ordered distinct $n$ points on $S_g$ as given in section $1$.
As above, fix a base-point
${\mathbf x}=(x_1, \cdots, x_n) \in {\rm {\rm Conf}}(S_g,n)$.
The fundamental group of the
configuration space
$\pi_1({\rm Conf}(S_g,n),{\mathbf x})$
is by definition the pure braid group
of $S_g$ with $n$ strands, and is denoted by $P_n(S_g)$.
The symmetric group $\Sigma_n$ acts freely on
${\rm Conf}(S_g,n)$ by permutation of components.
The fundamental group of the quotient space
${\rm Conf}(S_g,n)/\Sigma_n$ is the braid group of
$S_g$ with $n$ strands, and is denoted by $B_n(S_g)$.

Let us now assume $g>1$.
The fundamental group $\pi_1(S_g, x_0)$ with base-point $x_0$
is identified with a discrete subgroup
$G$ of $PSL(2, \mathbb R)$ acting freely on
$\mathbb H^2$. The $n$-fold product $G^n$ acts on the
orbit configuration space yielding a
covering space projection
$$
{\rm Conf}^{G}(\mathbb H^2,n) \to\ {\rm
Conf}(S_g,n)
$$
with the covering transformation
group $G^n$. In addition,
${\rm Conf}(S_g,n)$ is naturally a subspace of the $n$-fold
cartesian product $S_g^n$. Thus,
the inclusion map
$i:{\rm Conf}(S_g,n)
\to\ S_g^n$
induces a natural homomorphism
$$
i_* : P_n(S_g) \rightarrow \times_{j=1}^n \pi_1(S_g, x_j).
$$
Let
$P_n(S_g)^0$ denote the kernel of $i_*$.
By the above remarks, $P_n(S_g)^0$ is isomorphic to the fundamental
group of the orbit configuration space ${\rm Conf}^{G}(\mathbb H^2,n)$.
In the case $g=1$, we represent $S_g$ as an elliptic curve
${\mathbb C} / L$ with a lattice $L$
and we define $P_n(S_g)^0$ to be the fundamental group of the
orbit configuration space ${\rm Conf}^{L}({\mathbb C},n)$.
Thus, for any $g \geq 1$, we have an exact sequence
$$
1 \rightarrow P_n(S_g)^0 \rightarrow  P_n(S_g)
\rightarrow
\times_{j=1}^n \pi_1(S_g, x_j)
\rightarrow 1.
$$

Let us recall that any element of $P_n(S_g)$ is represented as
a collection of disjoint smooth paths
$\gamma = (\gamma_1, \cdots, \gamma_n)$ in $S_g \times [0, 1]$,
where the $i$-th string $\gamma_i$ runs
from the point $(x_i, 0)$ to some point $(x_j, 1)$.
Such notion of braids is extended to that of singular braids in
the following way.
We shall say that a collection of smooth paths
$\gamma = (\gamma_1, \cdots, \gamma_n)$ in $S_g \times [0, 1]$,
where the $i$-th string $\gamma_i$ runs monotonically in $t \in [0, 1]$
from the point $(x_i, 0)$ to some point $(x_j, 1)$,
is called a singular braid if the paths
$\gamma_1, \cdots, \gamma_n$
are allowed to intersect transversely, but only at finitely many
double points.
The set of isotopy classes of such singular braids is denoted by
${\mathcal S}B_n(S_g)$. Let us notice that
${\mathcal S}B_n(S_g)$ is not a group, but
has a structure of a monoid.
We denote by
${\mathbb Z}[B_n(S_g)]$ the group ring of $B_n(S_g)$.
We have a map
$$
\eta : {\mathcal S}B_n(S_g) \rightarrow
{\mathbb Z}[B_n(S_g)]
$$
defined in the following manner.
Let $\gamma$ be an element of ${\mathcal S}B_n(S_g)$ and we denote by
$p_1, \cdots, p_k$ its double points.
For $\epsilon_1=\pm 1, \cdots,
\epsilon_k = \pm 1$,
$\gamma_{\epsilon_1 \cdots \epsilon_k}$ stands for the braid
obtained by replacing each double point $p_i$
by a positive or a negative crossing according as
$\epsilon_i$ is $1$ or $-1$. For the above $\gamma \in {\mathcal S}B_n(S_g)$
we define $\eta$ by
\begin{equation}\label{eqn:singular braid}
\eta(\gamma) = \sum_{\epsilon_1=\pm 1, \cdots,
\epsilon_k = \pm 1} \epsilon_1 \cdots \epsilon_k
\ \gamma_{\epsilon_1 \cdots \epsilon_k}.
\end{equation}
Let ${\mathcal S}_kB_n(S_g)$ denote the set of isotopy classes of
singular braids with at least $k$ double points and we define
${\mathcal F}_k$ to be the ${\mathbb Z}$ submodule of
${\mathbb Z}[B_n(S_g)]$ generated by $\eta ({\mathcal S}_kB_n(S_g))$.
We have ${\mathcal F}_i {\mathcal F}_j \subset {\mathcal F}_{i+j}$
and it can be shown that each ${\mathcal F}_k$ is a two-sided ideal of
${\mathbb Z}[B_n(S_g)]$.
Thus, we obtain a decreasing filtration
$$
{\mathbb Z}[B_n(S_g)]=
{\mathcal F}_0 \supset {\mathcal F}_1 \supset
\cdots \supset {\mathcal F}_k \supset \cdots,
$$
which we shall call the Vassiliev filtration.
In a similar way, by putting ${\mathcal F}_k'={\mathcal F}_k \cap
{\mathbb Z}[P_n(S_g)]$, we obtain the Vassiliev filtration
$$
{\mathbb Z}[P_n(S_g)]=
{\mathcal F}_0' \supset {\mathcal F}_1' \supset
\cdots \supset {\mathcal F}_k' \supset \cdots,
$$
for the pure braid group of $S_g$.
A map $v: B_n(S_g) \rightarrow {\mathbb Z}$ is extended linearly
to a homomorphim of ${\mathbb Z}$-modules
$v : {\mathbb Z}[B_n(S_g)] \rightarrow {\mathbb Z}$.
We shall say that $v$ is a Vassiliev invariant
of order $k$
if $v$ vanishes on the ideal ${\mathcal F}_{k+1}$.
The set of Vassiliev invariants of order $k$ with values in
${\mathbb Z}$ is identified with
$$
\operatorname{Hom}_{\mathbb Z}
({\mathbb Z}[B_n(S_g)] / {\mathcal F}_{k+1} , {\mathbb Z})
$$
and we denote this set by
${\mathcal V}_k(B_n(S_g))$.
There is an increasing filtration
$$
{\mathcal V}_0(B_n(S_g)) \subset
{\mathcal V}_1(B_n(S_g)) \subset \cdots \subset
{\mathcal V}_k(B_n(S_g)) \subset \cdots
$$
A Vassiliev invariant of order $k$ for $P_n(S_g)$ is defined in
a similar way and the set of such invariants is identified with
$$
\operatorname{Hom}_{\mathbb Z}
({\mathbb Z}[P_n(S_g)] / {\mathcal F}_{k+1}' , {\mathbb Z}),
$$
which is denoted by
${\mathcal V}_k(P_n(S_g))$.

Our next object is to describe the weight systems for
Vassiliev invariants of braids on surfaces.
For this purpose we need the notion of horizontal chord diagrams on
surfaces.
Let $I_1 \sqcup \cdots \sqcup I_n$ be the disjoint union of
$n$ unit intervals. Fix a parametrization
$$
p_j : [0, 1] \rightarrow I_j
$$
for each $j$, $1\leq j \leq n$. A horizontal chord diagram on $n$ strands with
$k$ chords is a 1-dimensional complex
constructed in the following way. Fix
$t_1, \cdots, t_k \in [0, 1]$ such that $0 < t_1 < \cdots < t_k < 1$.
Let
$$
(i_1, j_1), (i_2, j_2), \cdots, (i_k, j_k)
$$
be pairs of distinct integers such that $1 \leq i_p \leq n$,
$1 \leq j_p \leq n$, $p=1, 2, \cdots, n$.
Consider $k$ copies of parametrized
unit intervals $C_1, C_2, \cdots, C_k$ and attach each
$C_{\nu}$ to $I_1 \sqcup \cdots \sqcup I_n$ in such a way that it starts at
$p_{i_{\nu}}(t_{\nu})$ and ends at $p_{j_{\nu}}(t_{\nu})$
for $1 \leq \nu \leq k$. In this way, we obtain a $1$-dimensional complex
with $n$ strands $I_1 \sqcup \cdots \sqcup I_n$ and
chords $C_1, \cdots, C_k$ attached to them.
Such a 1-dimensional complex is called a horizontal
chord diagram on $n$ strands with $k$ chords.
Here each strand and each chord are oriented by the
above parametrizations.
In what follows below, chord diagrams up to orientation preserving
homeomorphism
will be considered.

Let $C_n^k$ denote the set of horizontal chord diagrams
on $n$ strands with $k$ chords. For $\Gamma \in C_n^k$
consider a continuous map
$
f : \Gamma \rightarrow S_g
$
such that
\begin{equation}\label{eqn:basepoint}
f(p_i(0))=f(p_i(1))=x_i,  \  1 \leq i \leq n
\end{equation}
and denote by
$[f]$ its homotopy class.
Here, we consider the homotopy preserving the
condition \eqref{eqn:basepoint}.
Let $D_n^k(S_g)$ denote the free
${\mathbb Z}$-module
spanned by pairs
$(\Gamma, [f])$
for $\Gamma \in C_n^k$ and
$
f : \Gamma \rightarrow S_g,
$
a continuous map with the condition
\eqref{eqn:basepoint}.
The subspace of
$D_n^k(S_g)$ spanned by
$
f : \Gamma \rightarrow S_g,
$
such that each curve
$f(p_i(t)), 0 \leq t \leq 1$, is homotopic to
the point $\{ x_i \}$ is denoted by
$D_n^k(S_g)^0$.

Let ${\Gamma}_{ij}$ denote a horizontal chord diagram
on $n$ strands with one chord $C_{ij}$ defined by
the pair $(i, j)$, $i\neq j$, $1\leq i, j \leq n$.
Recall the base-point $x_0 \in S_g$ and
consider the fundamental group $\pi_1(S_g, x_0)$.
Next fix a path in $S_g$
connecting $x_0$ to $x_j$, and identify the set of
homotopy classes of paths from
$x_i$ to $x_j$ with
$\pi_1(S_g, x_0)$.
For $\gamma \in \pi_1(S_g, x_0)$ consider
$({\Gamma}_{ij}, [f]) \in D_n^k(S_g)^0$
such that $f(C_{ij})$ corresponds to
$\gamma \in \pi_1(S_g, x_0)$
by the above identification, and denote this
choice for $({\Gamma}_{ij}, [f])$
by $X_{i,j}^{\gamma}$.

Notice that the direct sum
$$
D_n(S_g)^0 = \bigoplus_{k \geq 0}D_n^k(S_g)^0
$$
has a structure of an algebra over ${\mathbb Z}$
where the product is defined by the concatenation
of chord diagrams. As an algebra
$D_n^k(S_g)^0$ is generated by $X_{i,j}^{\gamma}$, $1 \leq i \neq j \leq n$,
$\gamma \in \pi_1(S_g, x_0)$.
  From the definition it follows immediately that the
relation
$X_{i,j}^{\gamma}=X_{j,i}^{\gamma^{-1}}$
holds.
It turns out that $D_n(S_g)^0$ is a non-commutative associative algebra
freely generated by $X_{i,j}^{\gamma}$,
$1 \leq i<j \leq n$ over ${\mathbb Z}$.
The direct sum
$$
D_n(S_g) = \bigoplus_{k \geq 0}D_n^k(S_g)
$$
has a structure of an associative algebra as well. For
the subspace $D_n^0(S_g)$ spanned by the chord diagrams with
empty chord, there is a natural injection
$$
\iota_j :
\pi_1(S_g, x_j)
\rightarrow  D_n^0(S_g), \ 1 \leq j \leq n,
$$
which induces an isomorphism of
$\mathbb Z$ algebras
$$
{\mathbb Z}[G^n]
\cong
    D_n^0(S_g)
$$
with $G=\pi_1(S_g, x_0)$.
Here we fix an isomorphism
$\pi_1(S_g, x_0) \cong \pi_1(S_g, x_j)$
by means of a path from $x_0$ to $x_j$.
The algebra
$$
\Lambda_n=
{\mathbb Z}[G^n]
$$
acts on
$D_n(S_g)^0$ by the conjugation
$$
\Gamma \mapsto \iota_j(\mu) \ \ \Gamma \ \iota_j(\mu^{-1}),  \
\Gamma \in D_n(S_g)^0, \
\mu \in \pi_1(S_g, x_0).
$$

\begin{lem}\label{lem:conjugation}
With respect to the above action, the following holds.
\begin{align*}
&\iota_l(\mu) \ X_{i,j}^{\gamma} \
\iota_l(\mu^{-1}) = X_{i,j}^{\gamma} \ \ \text{for}  \ l\neq i,
j, \\
&\iota_i(\mu) \ X_{i,j}^{\gamma} \
\iota_i(\mu^{-1}) = X_{i,j}^{\mu \gamma},
   \\
&\iota_j(\mu) \ X_{i,j}^{\gamma} \
\iota_j(\mu^{-1}) = X_{i,j}^{\gamma \mu^{-1}}.
\end{align*}
\end{lem}

\begin{proof}
Consider the concatenation of the chord diagrams
$\iota_l(\mu)$ and $X_{i,j}^{\gamma}$,
where $\iota_l(\mu)$ is a chord diagram with
empty chord and the chord of $X_{i,j}^{\gamma}$
is considered to be the element $\gamma$
in $\pi_1(S_g, x_0)$. We construct a homotopy
of chord diagrams such that the initial point
of the above chord slides along the loop
$\mu \in \pi_1(S_g, x_0)$ in the negative
direction.
The resulting chord diagram
$X_{i,j}^{\mu \gamma} \, \iota_l(\mu)$
is homotopic to
$\iota_l(\mu) X_{i,j}^{\gamma}$.
This shows the second equality.
The other equalities are shown in
a similar way.
\end{proof}

The algebra $D_n(S_g)$ is considered to be the semidirect product
$$
D_n(S_g)^0 \rtimes \Lambda_n
$$
by the above action of $\Lambda_n$ on $D_n(S_g)^0$.
Write $X_{i,j}$ for $X_{i,j}^{e}$.
Let ${\mathcal I}$ be the two-sided ideal of
$D_n(S_g)$ generated by
\begin{align*}
&[X_{i,j}, X_{s,t}], \ \ {i, j, s, t} \ {\rm distinct}, \\
&[X_{i,j}, X_{j,s}+X_{i,s}],   \ \ {i, j, s} \ {\rm distinct}.
\end{align*}
Then we set
$$
{\mathcal A}_n(S_g) = D_n(S_g) / {\mathcal I},   \ \
{\mathcal A}_n(S_g)^0 = D_n(S_g)^0 / {\mathcal I} \cap D_n(S_g)^0.
$$
There is an action of $\Lambda_n$ on
${\mathcal A}_n(S_g)^0$ by conjugation and the semidirect
product
$${\mathcal A}_n(S_g)^0 \rtimes \Lambda_n$$ with respect to
this action is isomorphic to
${\mathcal A}_n(S_g)$.

\begin{lem}\label{lem:infinitesimal surface braid relations}
The algebra ${\mathcal A}_n(S_g)^0$ is generated by
$X_{i,j}^{\gamma}$, $1\leq i \neq j \leq n$,
$\gamma \in \pi_1(S_g, x_0)$, with relations:
\begin{align*}
&X_{i,j}^{\gamma}=X_{j,i}^{\gamma^{-1}} \\
&[X_{i,j}^{\gamma}, X_{s,t}^{\delta}]=0, \ \ {i, j, s, t} \ \text{distinct}, \\
&[X_{i,j}^{\gamma}, X_{j,s}^{\delta}+X_{i,s}^{\gamma\delta}]=0  \ \ {i, j, s} \
\text{distinct}.
\end{align*}
\end{lem}

\begin{proof}
The first two relations clearly hold from the definition.
Let us show that the the relation
$[X_{i,j}^{\gamma}, X_{j,s}^{\delta}+X_{i,s}^{\gamma\delta}]=0$ holds
in ${\mathcal A}_n(S_g)^0$
for distinct ${i, j, s}$.
  From Lemma \ref{lem:conjugation} we obtain
$$
[X_{i,j}^{\gamma}, X_{j,s}^{\delta}]
=\iota_i(\gamma) \iota_s(\delta^{-1})
[X_{i,j}, X_{j,s}]
\iota_i(\gamma^{-1}) \iota_s(\delta).
$$
Similarly, we have
$$
[X_{i,j}^{\gamma}, X_{i,s}^{\gamma\delta}]
=\iota_i(\gamma) \iota_s(\delta^{-1})
[X_{i,j}, X_{i,s}]
\iota_i(\gamma^{-1}) \iota_s(\delta).
$$
Hence the relation follows from
$[X_{i,j}, X_{j,s}+X_{i,s}]=0$.
Since $D_n(S_g)$ is the semidirect product
$D_n(S_g)^0 \rtimes \Lambda_n$ the ideal
${\mathcal I} \cap D_n(S_g)^0$ is generated by
the conjugation by $\Lambda_n$ of the
generators of ${\mathcal I}$.
This shows the desired statement.
\end{proof}

The algebras  ${\mathcal A}_n(S_g)^0$ and ${\mathcal A}_n(S_g)$
were introduced in \cite{GP} by the above generators and
relations. These algebras are considered to be graded algebras
by defining
$\deg X_{i,j}^{\gamma}=1$, $1 \leq i \neq j \leq n$,
$\gamma \in \pi_1(S_g, x_0)$,
and $\deg g=0$ for any $g \in \Lambda_n$.
There is a direct sum decomposition
$$
{\mathcal A}_n(S_g)^0 = \bigoplus_{k\geq 0}{\mathcal A}_n^k(S_g)^0, \ \
{\mathcal A}_n(S_g) = \bigoplus_{k\geq 0}{\mathcal A}_n^k(S_g),
$$
where ${\mathcal A}_n^k(S_g)^0$ and
${\mathcal A}_n^k(S_g)$ denote the degree $k$ parts.
Let $v$ be a Vassiliev invariant of order $k$ for $P_n(S_g)$.
Then there is a ${\mathbb Z}$-module homomorphism
$w(v) : D_n(S_g) \rightarrow {\mathbb Z}$
defined in the following way.
Let us take $\delta \in D_n^k(S_g)$ represented by
$f : \Gamma \rightarrow S_g$, where
$\Gamma$ is a horizontal chord diagram on $n$ strands with $k$ chords.
Contracting each chord on $S_g$, we may consider
$\delta$ as a projection diagram of a singular pure braid on $S_g$ with
$k$ transverse double points.
We set
\begin{equation}\label{eqn:weight}
w(v)(\delta) =
\sum_{\epsilon_1=\pm 1, \cdots,
\epsilon_k = \pm 1} \epsilon_1 \cdots \epsilon_k
\ v(\delta_{\epsilon_1 \cdots \epsilon_k}),
\end{equation}
where $\delta_{\epsilon_1 \cdots \epsilon_k}$ stands for the pure
braid on $S_g$ obtained from $\delta$ by replacing
each double point with positive or negative crossing
according as $\epsilon_j$ is $1$ or $-1$ for $1 \leq j \leq k$
as in
(\ref{eqn:singular braid}).
Since $v$ is an order $k$ invariant the above expression is well defined.
Thus, we obtain a $\mathbb Z$-module homomorphism
$$
w : {\mathcal V}_k(P_n(S_g)) \rightarrow
\operatorname{Hom}_{\mathbb Z}(D_n^k(S_g), {\mathbb Z}).
$$

\begin{lem}\label{lem:4T}
For $v \in  {\mathcal V}_k(P_n(S_g))$
the associated map
$w(v)$ vanishes on
the ideal $\mathcal I$
and induces a homomorphism
$$
w :{\mathcal V}_k(P_n(S_g)) \rightarrow
\operatorname{Hom}_{\mathbb Z}({\mathcal A}_n^k(S_g), {\mathbb Z})
$$
\end{lem}

\begin{proof}
Consider the $4$ chord diagrams
$\Gamma_j$, $1 \leq j \leq 4$, in $D_n^k(S_g)$ expressed as
\begin{align*}
&\Gamma_1=Y X_{i,j}X_{j,s} Z, \ \Gamma_1=Y X_{j,s}X_{i,j} Z \\
&\Gamma_3=Y X_{i,j}X_{i,s} Z, \ \Gamma_4=Y X_{i,s}X_{i,j} Z
\end{align*}
with some fixed $Y, Z \in D_n(S_g)$.
Contracting the chords on $S_g$, we represent
$\Gamma_j$, $1 \leq j \leq 4$, as singular pure braids on $S_g$.
Replacing the double points by positive or negative crossings
and applying the definition \ref{eqn:weight}, we see that
the relation
$$
w(v)(\Gamma_1)-
w(v)(\Gamma_2)+
w(v)(\Gamma_3)-
w(v)(\Gamma_4)=0
$$
holds. This completes the proof.
\end{proof}

The above $w(v)$ is called the weight system of the
Vassiliev invariant $v$.
To define the weight system of
a Vassiliev invariant for
$B_n(S_g)$, we need to extend the algebra
${\mathcal A}_n(S_g)$ by the group algebra of the
symmetric group $\Sigma_n$ in the following way.
We define $\widetilde{\Lambda}_n$ to be the group algebra
of the semidirect product
$$
G^n  \rtimes \Sigma_n
$$
where the symmetric group acts as permutations.
We define $\widetilde{{\mathcal A}}_n(S_g)$ as the
semidirect product
$$
{\mathcal A}_n(S_g)^0 \rtimes \widetilde{\Lambda}_n.
$$
where the action of $\Sigma_n$ on ${\mathcal A}_n(S_g)^0$
is defined by
$$
\sigma X_{i,j}^{\gamma} \sigma^{-1}
= X_{\sigma(i),\sigma(j)}^{\gamma},  \ \sigma \in \Sigma_n
$$
and the action of $G^n$
is the conjugation action defined before.
The algebra $\widetilde{{\mathcal A}}_n(S_g)$
has a structure of a graded algebra as an extension of
the graded algebra
${\mathcal A}_n(S_g)$ by defining
$\deg g =0$ for any $g \in \widetilde{\Lambda}_n$.
We denote the degree $k$ part by
$\widetilde{{\mathcal A}}^k_n(S_g)$.
The weight system for $B_n(S_g)$ is defined by a natural
${\mathbb Z}$-module homomorphism
$$
w :{\mathcal V}_k(B_n(S_g)) \rightarrow
\operatorname{Hom}_{\mathbb Z}
(\widetilde{{\mathcal A}}_n^k(S_g), {\mathbb Z}).
$$

Let us consider the graded algebras
$\operatorname{gr}{\mathbb Z}[B_n(S_g)]$ and
$\operatorname{gr}{\mathbb Z}[P_n(S_g)]$
associated with the Vassiliev filtration.
Namely, we set
$$
\operatorname{gr}{\mathbb Z}[B_n(S_g)] =
\bigoplus_{k\geq 0} {\mathcal F}_{k} / {\mathcal F}_{k+1}, \ \
\operatorname{gr}{\mathbb Z}[P_n(S_g)] =
\bigoplus_{k\geq 0} {\mathcal F}_{k}' / {\mathcal F}_{k+1}'.
$$
The following theorem was shown by
J.~Gonz\'alez-Meneses and L. Paris.
\begin{thm}[\cite{GP}]\label{thm:GP}
There exists a homomorphism of ${\mathbb Z}$-modules
$u : {\mathbb Z}[B_n(S_g)] \rightarrow
\widetilde{{\mathcal A}}_n(S_g)$ such that the associated graded map
gives an isomorphism of graded ${\mathbb Z}$ algebras
$$
\operatorname{gr}
{\mathbb Z}[B_n(S_g)] \cong \widetilde{{\mathcal A}}_n(S_g).
$$
\end{thm}

The restriction of the above homomorpshism $u$ to
${\mathbb Z}[P_n(S_g)]$ gives an isomorphism of graded ${\mathbb Z}$-algebras
$$
\operatorname{gr}
{\mathbb Z}[P_n(S_g)] \cong {\mathcal A}_n(S_g).
$$
In particular, there is an isomorphism of ${\mathbb Z}$-modules
$$
{\mathcal F}_k' / {\mathcal F}_{k+1}'
\cong
{\mathcal A}_n^k(S_g).
$$

Now the statement of Theorem \ref{thm:chord diagrams}
follows from Theorem \ref{thm:Lie algebras} together with
Lemma \ref{lem:infinitesimal surface braid relations}.
It can also be shown directly using an argument due to
J.~Gonz\'alez-Meneses and L. Paris in the proof of
Theorem \ref{thm:GP} in the following way.
The restriction of the Vassiliev filtration for
${\mathbb Z}[P_n(S_g)]$ on
${\mathbb Z}[P_n(S_g)^0]$
coincides with the filtration given by the powers of the
augmentation ideal $I$ of ${\mathbb Z}[P_n(S_g)^0]$.
Hence, the associated graded algebra
$$
{\rm gr}_I {\mathbb Z}[P_n(S_g)^0] = \bigoplus_{k\geq 0}I_k/I_{k+1}
$$
is isomorphic to
the graded algebra ${\mathcal A}_n(S_g)^0$.
It can be shown that the successive quotients of the descending
central series $\Gamma^m(P_n(S_g)^0) / \Gamma_{m+1}(P_n(S_g)^0)$, $m=0, 1,
\cdots$,
is a free ${\mathbb Z}$-module.
Therefore, using a theorem due to Quillen \cite{Q}, we obtain
the isomorphism of graded algebras
$U[E_0^*(P_n(S_g)^0)] \cong {\rm gr}_I {\mathbb Z}[P_n(S_g)^0]$.
This shows the isomorphism of graded algebras
$U[E_0^*(P_n(S_g)^0)] \cong {\mathcal A}_n(S_g)^0$.

The exceptional case of $g=1$ was considered in articles \cite{K3,CX}.
This case arises by considering horizontal chord diagrams
for elliptic curves \cite{K3}. Analogous information associated
to the orbit configuration space ${\rm Conf}^L(\mathbb C,n)$
where $L$ denotes a fixed lattice in the complex plane
$\mathbb C$ such that $\mathbb C / L$ is a fixed elliptic curve
is given in \cite{CX}.

\section{The Lie algebra relations in Theorem \ref{thm:Lie algebras}}

The purpose of this section is to define the elements
$B^{\tau}_{i,j}$, and to derive the relations among these
elements as stated in Theorem $1$($3$). To
analyze this structure, it is convenient
to construct explicit cycles in $ {\rm Conf}^{G}(\mathbb H^2 \times \mathbb
C^q,n) $
as follows.

Choose points $p_1,p_2, \cdots, p_n$ in distinct orbits of $G$ in
$\mathbb H^2 \times \mathbb C^q$. Let $\sigma$, and $\tau$ denote
elements of $G$. Given fixed $i$, and $j$ fixed with
$ 1 \leq j < i \leq n $, and $\sigma$ in $G$, choose a closed neighborhood $U$
homeomorphic to a closed disk with the following properties.

\begin{enumerate}
    \item $U$ contains the interior points $p_j$,
    together with the points $\sigma(p_j)+ z/ \lambda$ for $\lambda$
    a fixed real scalar for all $z$ in $\mathbb H^2 \times \mathbb C^q$
    of norm $1$ where the hyperbolic metric is used for $\mathbb H^2$,
    \item $U$ intersects each orbit exactly once.
    \item $U$ does not contain the $p_t$ for all $t$ with $t \neq j$.
\end{enumerate}

Define maps $$A^{\sigma}_{i,j}: S^{2q+1} \to\
{\rm Conf}^{G}(\mathbb H^2 \times \mathbb C^q,n)$$ by the formula
$$A^{\sigma}_{i,j}(z) =
(w_1,w_2, \cdots,w_k)$$ where

\begin{enumerate}
\item $w_t = p_t$ if $t \neq i $, and
\item $w_i = \sigma(p_j) + z/ \lambda $.
\end{enumerate}

Fix a choice of fundamental cycle $\iota$ for $S^{2q+1}$ for what follows
below. The image in homology ${{A^{\sigma}}_{i,j}}_*(\iota)$
will be denoted by the ``name" $A^{\sigma}_{i,j}$
in $H_{2q+1}{\rm Conf}^{G}(\mathbb H^2 \times \mathbb C^q,n)$.
In addition, the ``name" $B^{\sigma}_{i,j}$ is also used ambiguously
for

\begin{enumerate}
    \item the analogous element in $\pi_1({\rm Conf}^{G}(\mathbb H^2,n),\mathbf
    x)$ denote , and
    \item the image of the transgression of the spherical class
$A^{\sigma}_{i,j}$ in
    $H_{2q}\Omega {\rm Conf}^{G}(\mathbb H^2 \times \mathbb C^q,n)$.
\end{enumerate}

Relations among the elements $A^{\sigma}_{i,j}$ are obtained as follows.
The first relation is given by the next lemma.

\begin{lem}\label{lem:Lie antisymmetry}
The relation $$\tau_{(i,j)}A^{{\sigma} }_{i,j} = A^{{\sigma}^{-1} }_{i,j}$$
holds in
\begin{align*}
&H_{2q+1}{\rm Conf}^{G}(\mathbb H^2 \times \mathbb C^q,n), and \\
& E_0^*(\pi_1({\rm Conf}^{G}(\mathbb H^2, n),\mathbf x))\
\end{align*} where $\tau_{(i,j)}$ denotes the transposition
which interchanges $i$, and $j$ in the symmetric group $\Sigma_n$.

Furthermore, the action of an element $\gamma$ in the symmetric group on
the element
$A^{{\sigma} }_{i,j}$ is specified by the formula

\[
\gamma(A^{{\sigma} }_{i,j}) =
\begin{cases}
   A^{{\sigma} }_{\gamma(i),\gamma(j)} & \text{if $\gamma(i)> \gamma(j)$,}\\
A^{{\sigma}^{-1} }_{\gamma(j),\gamma(i)} & \text{if $\gamma(i)< \gamma(j)$.}
\end{cases}
\]

\end{lem}

\begin{proof}

First recall the standard homotopy for the classical configuration
space with $G = \{identity\}$. Fix points $p_m = 4m{e_1}$, in $\mathbb C^q$
for $1 \leq m \leq n$ for which $e_1$ is a unit vector. Define
$$h:[0,1] \times S^{2q-1} \to\ {\rm Conf}( \mathbb C^q,n)$$ by the
formula $$h(t,z) = (x_1,x_2, \cdots,x_n)$$ for which
\begin{enumerate}
      \item  $x_s = p_s$ if $s \neq i,j$,
      \item $x_j = p_j - tz$,
      \item $x_i = p_i + (1-t)z$.
\end{enumerate}

Notice that the point $(x_1,x_2, \cdots,x_n)$ is in the
configuration space, and that the map $h$ is continuous.
In particular, each $x_s$ is a function of $t$ with
$x_j(0) = p_j$ with
\begin{enumerate}
    \item $x_j(0) = p_j$,
    \item $x_j(1) = p_j -z$,
    \item $x_i(0) = p_i +z$, and
    \item $x_i(1) = p_i$.
\end{enumerate} Thus the maps
\begin{enumerate}
    \item $\tau(A_{i,j}): S^{2q-1} \to\ {\rm Conf}( \mathbb C^q,n)$,
    and
    \item $A_{i,j}: S^{2q-1} \to\ {\rm Conf}( \mathbb C^q,n)$
\end{enumerate} are homotopic.

An analogous homotopy is constructed next for $A^{\sigma}_{i,j}$
as follows. The small disk $U$ is regarded as
the standard unit disk in $\mathbb R^{2q+2}$. Define
$$H:[0,1] \times S^{2q+1} \to\ {\rm Conf}^{G}(\mathbb H^2 \times \mathbb
C^q,n)$$
by the formula

$$H(t,z) = (y_1,y_2, \cdots,y_n)$$
for which

\begin{enumerate}
      \item  $y_s = p_s$ if $s \neq i,j$,
      \item $y_j = p_j - tz/\lambda$, and
      \item $y_i = \sigma(p_j) + (1-t)z/\lambda$.
\end{enumerate} Hence $H$ is a homotopy of $A^{\sigma}_{i,j}$
(a homotopy which is non-base-point preserving). Notice that
$(\tau_{(i,j)}\circ H)(t,z) = (u_1,u_2,\cdots, u_k)$
where

\begin{enumerate}
      \item $u_s = p_s$ if $s \neq i,j$,
      \item $u_j = \sigma(p_j) + (1-t)z/\lambda$, and
      \item $u_i = p_j - tz/\lambda$.
\end{enumerate} Thus $\tau_{(i,j)}(A^{\sigma}_{i,j})$ is homotopic to the map
which sends $z$ to $(u_1,u_2,\cdots, u_n)$ where

\begin{enumerate}
      \item $u_s = p_s$ if $s \neq i,j$,
      \item $u_j = \sigma(p_j)$, and
      \item $u_i = p_j - z/\lambda$.
\end{enumerate} which is given by $A^{\sigma^{-1}}_{i,j}$
up to (a non-base-point preserving ) homotopy.

A similar computation applies to $\gamma(A^{{\sigma} }_{i,j})$, and
is deleted. This suffices.
\end{proof}

Analogous maps, and homotopies are considered next.
Fix integers $1 \leq j < s < i \leq n $, and define maps
$$D(\sigma, \tau, i,s,j): S^{2q+1} \times S^{2q+1} \to\
{\rm Conf}^{G}(\mathbb H^2 \times \mathbb C^q,n), $$ and
$$G(\sigma, \tau,i,s,j): S^{2q+1} \times S^{2q+1} \to\
{\rm Conf}^{G}(\mathbb H^2 \times \mathbb C^q,n)$$
by the following formulas:
\begin{enumerate}
\item $$D(\sigma, \tau, i,s,j)(u,v) = (w_1,w_2, \cdots, w_n)$$ where
      \begin{enumerate}
      \item  $w_t = p_t$ if $t \neq i,s,j $,
      \item  $w_j = p_j $,
      \item  $w_s = \sigma (p_j)+ u/8 \lambda $,
      and
      \item $w_i = \tau(p_j) + v/16 \lambda  $.
      \end{enumerate}

\item
   $$G(\sigma, \tau, i,s,j)(u,v) = (w_1,w_2, \cdots, w_n)$$ where
      \begin{enumerate}
      \item $w_t = p_t$ if $t \neq i,s,j $,
      \item $w_j = p_j$,
      \item $w_s = \sigma (p_j) + v/16 \lambda $, and
      \item $w_i = \tau (p_j) + u/8 \lambda $.
      \end{enumerate}
\end{enumerate}

Write $a \otimes 1$, and $1 \otimes b$ for
the respective generators in $H_{2q+1}(S^{2q+1} \times S^{2q+1})$
associated to  the individual fundamental cycles of each axis.
Consider the induced map in homology for
$$D(\sigma, \tau, i,s,j)_*: H_{2q+1}(S^{2q+1} \times S^{2q+1}) \to\
H_{2q+1}{\rm Conf}^{G}(\mathbb H^2 \times \mathbb C^q,n).$$  Notice that
$H_{2q+1}{\rm Conf}^{G}(\mathbb H^2 \times \mathbb C^q,n)$ is
free abelian with basis
$$\{A^{\sigma}_{i,j} \ | \ \sigma \in
G, 1 \leq j < i \leq n\}$$
and where the projection to the
$i$, and $j$ coordinates $$p_{i,j}:{\rm Conf}^{G}(\mathbb H^2 \times
\mathbb C^q,n) \to\
{\rm Conf}^{G}(\mathbb H^2 \times \mathbb C^q,2)$$ is non-trivial in
homology when restricted to
$A^{\sigma}_{i,j}$. Thus projecting gives the following formulas:

\begin{enumerate}
            \item $p_{i,j} D(\sigma, \tau, i,s,j)(u,v) =
            (p_j,\tau (p_j) + v/16 \lambda)$ which is homotopic
            to $A^{\tau}_{i,j}(v)$,

            \item $p_{i,s}\circ D(\sigma, \tau, i,s,j)(u,v) = (\sigma (p_j)
+ u/8
            \lambda, \tau(p_j) + v/16 \lambda )$ which is homotopic to
            $A^{\tau \sigma^{-1}}_{i,s}(u)$,

            \item $p_{s,j} D(\sigma, \tau, i,s,j)(u,v) = (p_j,\sigma (p_j) + u/8
            \lambda)$ which is homotopic
            to $A^{\sigma}_{s,j}(u)$, and

            \item $p_{ \alpha,\beta } D(\sigma, \tau, i,s,j)(u,v)$ is null
otherwise.
\end{enumerate}

\begin{lem}\label{lem:D relations}
The following formulas are satisfied:
\begin{enumerate}
            \item ${D(\sigma, \tau, i,s,j)}_*(a \otimes 1) =
            A^{\tau \sigma^{-1} }_{i,s}+ A^{\sigma}_{s,j}$,
            \item $D(\sigma, \tau, i,s,j)_*(1 \otimes b) = A^{\tau}_{i,j}$, and
\item $[B^{\tau \sigma^{-1} }_{i,s}+ B^{\sigma}_{s,j}, B^{\tau}_{i,j}] = 0$ in
$H_*(\Omega {\rm Conf}^{G}(\mathbb H^2 \times \C^q,n))$, and
in $E_0^*(\pi_1 {\rm Conf}^{G}(\mathbb H^2,n), \mathbf x)$.
\end{enumerate}
\end{lem}

\begin{proof}
Notice that formula $1$ follows from the above remarks.
Furthermore, in case $n > 0$, the elements $B^{\tau}_{i,j}$
are defined to be the adjoint of the elements $A^{\tau}_{i,j}$,
and thus the relation $2$ holds in homology.

In addition, formula $2$ is satisfied on the level of the
associated graded as the homotopy $D$ gives that
$[B^{\tau \sigma^{-1} }_{i,s}\cdot B^{\sigma}_{s,j},B^{\tau}_{i,j}]$ lies in
$\Gamma^2(\pi_1 {\rm Conf}^{G}(\mathbb H^2,n), \mathbf x)$, the second stage
of the descending central series. The lemma
follows.
\end{proof}

Similar formulas apply to $G(\sigma, \tau, i,s,j)(u,v)$.

\begin{enumerate}
            \item $p_{i,j} G(\sigma, \tau, i,s,j)(u,v) =
            (p_j,\tau (p_j) + u/8 \lambda )$ which is homotopic
            to $A^{\tau }_{i,j}(u)$,

            \item $p_{i,s} G(\sigma, \tau, i,s,j)(u,v) =
            (\sigma (p_j) + v/16\lambda, \tau (p_j) + u/8 \lambda )$
            which is homotopic to
            $A^{\tau \sigma^{-1}}_{i,s}(u)$,

            \item $p_{s,j} G(\sigma, \tau, i,s,j)(u,v) = (p_j,\sigma (p_j) +
            v/16 \lambda)$ which is homotopic
            to $A^{\sigma}_{s,j}(v)$, and

            \item $p_{ \alpha,\beta } G(\sigma, \tau, i,s,j)(u,v)$ is null
otherwise.

\end{enumerate}

In addition, the following formulas are satisfied:
\begin{enumerate}
            \item ${G(\sigma, \tau, i,s,j)}_*(a \otimes 1) =
            A^{\tau }_{i,j} + A^{\tau \sigma^{-1} }_{i,s}$,

            \item $G(\sigma, \tau, i,s,j)_*(1 \otimes b) =
            A^{\sigma}_{s,j}$, and

            \item $[B^{\sigma}_{s,j},B^{\tau }_{i,j} + B^{\tau
\sigma^{-1}}_{i,s}] = 0$.
\end{enumerate}\vskip .2in

The statement and proof of the next lemma are analogous.

\begin{lem}\label{G:relations}
The following formulas are satisfied:
\begin{enumerate}
            \item ${G(\sigma, \tau, i,s,j)}_*(a \otimes 1) =
            A^{\tau }_{i,j} + A^{\tau \sigma^{-1} }_{i,s}$,

            \item $G(\sigma, \tau, i,s,j)_*(1 \otimes b) =
            A^{\sigma}_{s,j}$, and

            \item $[B^{\sigma}_{s,j},B^{\tau }_{i,j} + B^{\tau
\sigma^{-1}}_{i,s}] = 0$ in
            $H_*(\Omega {\rm Conf}^{G}(\mathbb H^2 \times \C^q,n))$, and
            in $E_0^*(\pi_1 {\rm Conf}^{G}(\mathbb H^2,n), \mathbf x)$.
\end{enumerate}
\end{lem}

Define $$H(\sigma, \tau,i,j,s,t): S^{2q+1} \times S^{2q+1} \to\
{\rm Conf}^{G}(\mathbb H^2 \times \mathbb C^q,n)$$
by the following formulas where $ 1 \leq j < i <s < t$: $$H(\sigma,
\tau, i,s,j)(u,v) = (w_1,w_2, \cdots, w_n)$$ where
      \begin{enumerate}
      \item  $w_m = p_m$ if $m \neq i,j,s,t $,
      \item  $w_j = p_j $,
      \item  $w_i = \sigma (p_j)+ u/8 \lambda $,
      \item  $w_s = p_s $, and
      \item  $w_t = \tau (p_s)+ v/8 \lambda $.
      \end{enumerate}

Thus the following formulas are satisfied:
\begin{enumerate}
            \item ${H(\sigma, \tau, i,j,s,t)}_*(a \otimes 1) =
            A^{\sigma}_{i,j} $,

            \item ${H(\sigma, \tau, i,j,s,t)}_*(1 \otimes b) =
            A^{\tau }_{t,s} $,and

            \item $[B^{\sigma}_{i,j},B^{\tau }_{s,t}] = 0$
            in $H_*(\Omega {\rm Conf}^{G}(\mathbb
            H^2 \times \C^q,n) )$ if $ 1 \leq j < i <s < t$.

            \item The other cases for which $\{i,j\} \cap \{s,t\} = \phi$
            are omitted as the details are similar.
\end{enumerate}

\begin{lem}\label{H:relations}
The following formulas are satisfied:
\begin{enumerate}
   \item ${H(\sigma, \tau, i,j,s,t)}_*(a \otimes 1) =
            A^{\sigma}_{i,j} $,

            \item ${H(\sigma, \tau, i,j,s,t)}_*(1 \otimes b) =
            A^{\tau }_{t,s} $, and

            \item $[B^{\sigma}_{i,j},B^{\tau }_{s,t}] = 0$
            in $H_*(\Omega {\rm Conf}^{G}(\mathbb
            H^2 \times \C^q,n) )$ in
            $H_*(\Omega {\rm Conf}^{G}(\mathbb H^2 \times \C^q,n)
            )$, and in $E_0^*(\pi_1 {\rm Conf}^{G}(\mathbb
            H^2,n), \mathbf x)$ if $\{i,j\} \cap \{s,t\} = \phi$.
\end{enumerate}
\end{lem}

The proof of this lemma is similar to the previous two, and is
omitted. The relations stated in Theorem \ref{thm:Lie algebras} follow.

\section{ Proof of Theorem \ref{thm:Lie algebras} }

Recall the following analogue of results due to Fadell and Neuwirth \cite{FN}
for ``usual" configuration spaces, and which apply to orbit configuration
spaces \cite{X}.

\begin{lem}\label{L:fibrations}
Let $M$ be a manifold with a free, and properly discontinuous action
of a group $G$. Then for $j \leq n$, the projection
$p: {\rm Conf}^{G}(M,n) \to\ {\rm Conf}^{G}(M, j)$ onto the first $j$
coordinates  is a locally trivial bundle, with fibre
${\rm Conf}^{G}(M - Q_{j}^{G}, n-j)$.
\end{lem}

It does not suffice to use free actions in the above lemma.
The extra hypothesis that the action be properly discontinuous
is required for the case when $G$ is not finite. Notice
that the hypotheses are trivially satisfied when
$G$ is a finite group acting freely on a Hausdorff space $M$.

The next lemma guarantees cross-sections for these fibrations
in favorable cases.

\begin{lem}\label{L:sections}
Let $M$ be a parallelizable manifold with a free, and properly
discontinuous action
of a group $G$ such that the quotient map  $M \to\ M/ G$ is the projection for
a covering space. Then for $j < n$, the projection
$p: {\rm Conf}^{G}(M,n) \to\ {\rm Conf}^{G}(M, j)$ admits a
cross-section. Thus if $M$ is simply-connected of dimension
at least $3$, then there are homotopy equivalences
$$\Omega ({\rm Conf}^{G}(M,j)) \times \Omega({\rm Conf}^{G}(M - Q_{j}^{G},
n-j))
\to\ \Omega({\rm Conf}^{G}(M,n)).$$ Thus there is a homotopy
equivalence $$\times_{0 \leq j \leq n-1 } \Omega(M - Q_{j}^{G})
\to\ \Omega( {\rm Conf}^{G}(M,n)).$$
\end{lem}

\begin{proof}
The exponential map with source the tangent bundle of $M$, $exp: \tau(M) \to\
M$, gives a local homeomorphism from $\mathbb R^m$ to an
open set in $M$ (where $M$ is a manifold of dimension $m$ ).
Choose $n$ elements $y_1,y_2, \cdots, y_n$ in $\mathbb R^m$
in $M$ with $y_1 = 0$, such that $x_1,x_2, \cdots, x_n$ for $x_i = exp(y_i)$
lie in $n$ distinct orbits as the action of $G$ is properly discontinuous.
A cross-section is defined by by $\sigma(x)=
(x_1,x_2,\cdots, x_n)$ where $x_i = exp(y_i)$.

Notice that there are
fibrations $${\rm Conf}^{G}(\mathbb H^2 \times \mathbb C^q,n) \to\
{\rm Conf}^{G} (\mathbb H^2 \times \mathbb C^q,n-1)$$ with both
(1) cross-sections, and (2) fibres given by
$$\mathbb H^2 \times \mathbb C^q - Q_{n-1}^{G}.$$ Since
any multiplicative fibration with section is homotopy equivalent
to a product, there are homotopy equivalances
$$\prod_{0 \leq i \leq n-1} \Omega (\mathbb H^2
\times \mathbb C^q - Q_{i}^{G}) \to\
\Omega({\rm Conf}^{G}(\mathbb H^2 \times \mathbb C^q,n)$$ for $ q \geq 1 $.
The lemma follows.
\end{proof}

The previous lemma gives the proof of Theorem \ref{thm:Lie algebras} in
case $q > 0$.
That is the maps $ \Omega (\mathbb H^2 \times \mathbb C^q - Q_{i}^{G}) \to\
\Omega({\rm Conf}^{G}(\mathbb H^2 \times \mathbb C^q,n)$ for $ q \geq 1 $
induce multiplicative maps on the level of homology groups which
are split monomorphisms. The only relations satisfied are those
given already as the homology of $\Omega (\mathbb H^2
\times \mathbb C^q - Q_{i}^{G})$ is a tensor algebra.

It suffices to prove the theorem in case $q = 0$.
The proof in this case follows at once from Theorem \ref{thm:the case H}
which is proven in the next section.

\section{Proof of Theorem \ref{thm:the case H}}

Fix a base-point $y_n$ in ${\rm Conf}^{G}(\mathbb H^2,n)$.
The direct sum decomposition of Lie algebras given in
Theorem \ref{thm:the case H} follows from
the split, short exact sequence of groups with
trivial local coefficients given by
$$1 \to\ \pi_1(\mathbb H^2-Q_{n-1}^{G}, x) \to\
\pi_1({\rm Conf}^{G}(\mathbb H^2,n), \mathbf y_n)\to\
\pi_1({\rm Conf}^{G}(\mathbb H^2,n-1), \mathbf y_{n-1}) \to\ 1.$$
A split, short exact sequence of Lie algebras as
restated in the next lemma is given in \cite{K1, FR, Xico}.

\begin{lem}\label{splitting}
Let $$1 \to\ A \to\ B \to\ C \to\ 1$$  be a split short exact sequence
of groups for which the conjugation of $C$ on $H_1(A)$ is trivial.
Then there is a split short exact sequence of Lie algebras
$$0 \to\ E_0^{i*}(A) \to\ E_0^{i*}(B)  \to\ E_0^{i*}(C)  \to\ 0.$$
\end{lem}

Hence the proof of Theorem \ref{thm:the case H} rests on the proof
of the triviality of the local coefficient system for the fibration
$$ \mathbb H^2-Q_{n-1}^{G}, \to\ {\rm Conf}^{G}(\mathbb H^2,n) \to\
{\rm Conf}^{G}(\mathbb H^2,n-1),$$ a proof analogous to that given in \cite
{CX},
and in Lemma \ref{Local coefficients}.  Since the
homology of the fibre is concentrated in degree one, and the
fibration has a cross-section, the Serre spectral sequence
for these fibrations collapses. The next proposition gives
the additive structure for the integral homology
of ${\rm Conf}^G(\mathbb H^2, n)$ as stated in
Theorem \ref{thm:the case H}.

\begin{prop}\label{prop:additive homology}
The integral homology of ${\rm Conf}^G(\mathbb H^2, n)$ is additively given by
$$ H_* {\rm Conf}^G(\mathbb H^2, n) \cong  H_* (C_1) \otimes H_* (C_2) \otimes
\dots \otimes H_* (C_{n-1})$$ where $C_i$ is the infinite bouquet of circles
$\bigvee_{|Q_i^{G}|} S^1$.
\end{prop}

Notice that the homology of ${\rm Conf}^G(\mathbb H^2, n)$
is not of finite type. Thus the double dualization does not give an isomorphism
of vector spaces from the homology to the double dual of the homology.

\begin{prop}\label{L:integral homology}
The integral cohomology of  ${\rm Conf}^G(\mathbb H^2, n)$
is given additively by an isomorphism
$$ H^*{\rm Conf}^G(\mathbb H^2, n) \to\ H^* (C_1) \otimes \dots \otimes
H^* (C_{n-1}).$$
In addition, for every $i=2,\dots n$, there are choices of
classes in $H^*{\rm Conf}^G(\mathbb H^2, n)$ given by
$A^{\sigma_1}_{i,1},  A_{i,2}^{\sigma_2}, \dots, A^{\sigma_{i-1}}_{i,i-1}$
for $\sigma_p \in G$,
which are in a 1-1 correspondence with
the generators in $H_1(C_{i-1})$,
satisfying the following relations:

\begin{enumerate}
\item $A_{i,j}^{\mu} \cdot A_{i,j}^{\nu}  = 0$  for all
$\mu$, and  $\nu$  in $G$.
\item $A_{i,t}^{\mu} \cdot A_{i,j}^{\nu} =
{A}^{\mu {\nu}^{-1}}_{j, t} \cdot ({A}^{\nu}_{i,j} - {A}^{\mu}_{i,t})$ if
$ 1\leq t <j < i \leq k $ with   $\mu$ and, $\nu$
in $G$.
\end{enumerate}
\end{prop}

\begin{lem}\label{L:fundamental}
Let $p : E \rightarrow B$ be a locally trivial bundle with
$B$ path-connected. Let $x, y \in B$ and let $F_x$ and  $F_y$
be the corresponding fibers. Then $\pi_1(B, x)$ acts trivially
on $H_*(F_x)$ if and only if
$\pi_1(B, y)$ acts trivially on $H_*(F_y)$.
\end{lem}

\begin{lem}\label{Local coefficients}
The fibration
$$ \pi : {\rm Conf}^{G}( \mathbb H^2,n) \to\
{\rm Conf}^{G}( \mathbb H^2,n-1)$$
has trivial local coefficients.
\end{lem}

\begin{proof}

The construction at the beginning of section $3$
for points $p_1,p_2, \cdots, p_k$ in distinct orbits of $G$ in
$\mathbb H^2$ together with a closed neighborhood $U$
homeomorphic to a closed disk satisfies the following properties.

\begin{enumerate}
    \item $U$ contains the interior points $p_j$,
    together with the points $\sigma(p_j)+ z/ \lambda$ for $\lambda$
    a fixed real scalar for all $z$ in $\mathbb H^2 $
    of norm $1$ where the hyperbolic metric is used for $\mathbb H^2$,
    \item $U$ intersects each orbit exactly once.
    \item $U$ does not contain the $p_t$ for all $t$ with $t \neq j$.
\end{enumerate}

Next fix one value of $z$, say $\alpha$. Consider the Dehn twist
which ``rotates" $p_j$ to $\sigma(p_j) + \alpha/ \lambda $, and fixes
the boundary of $U$ pointwise. This disk, and Dehn twist is used to describe
the local coefficient system for certain fibrations below. This
Dehn twist is isotopic to the identity by an isotopy which fixes
the complement of $U$ pointwise.

The proof of the lemma is by downward induction on $r$ and is analogous to that
given in \cite{CX} in the special case of
${\rm Conf}^{L}( \mathbb C ,n)$ for the standard integral lattice
$L$ in $ \mathbb C$. Namely, the isotopy in the previous
paragraph gives a homeomorphism of $\mathbb H^2  - Q_i^{G}$ for which
$Q_i^{G}= \amalg_{1 \leq j \leq i} G \cdot p_j$. A direct
inspection gives that the effect of this homeomorphism on the level
of $H_1(\mathbb H^2 - Q_i^{G})$ is the identity,
but is not the identity on $\pi_1(\mathbb H^2 - Q_i^{G},x)$.
\end{proof}

\section{Graded Poisson algebra structures}

There is a natural structure of graded Poisson algebra for the homology
of an iterated loop space. It is the purpose of this section to
describe this structure for the case of $\Omega^k {\rm Conf}^{G}(\mathbb
H^2 \times
\mathbb C^q,n)$ for $k > 1$.

It is standard that a graded associative algebra $A$ inherits the structure
of a
graded Lie algebra with Lie bracket given by the commutator $$[a,b] = a\cdot b
- (-1)^{|a|\cdot|b|} b\cdot a$$ for an element $a$ of degree $|a|$, and
an element $b$ of degree $|b|$ in the graded associative algebra $A$. Thus the
homology of a $1$-fold loop space is naturally a graded Lie
algebra.

If $k>1$, the homology of a $k$-fold loop space, $\Omega^k X$,
admits the structure of a graded Poisson algebra \cite{CLM}, pages $215-217$.
Namely, there is a bilinear map, the Browder operation, given by
\[
\lambda_{k-1}: H_i(\Omega^k X) \otimes H_j(\Omega^k X) \to
H_{i+j+k-1} (\Omega^k X)
\] which satisfies the following axioms for a graded Poisson
algebra for which the degree of an element $x$ is written $|x|$:

\begin{enumerate}
    \item (Jacobi identity)
$$\alpha \lambda_{k-1}[a,\lambda_{k-1}[b,c]] + \beta \lambda_{k-1}[b,\lambda_{k-1}[c,a]] +
\gamma \lambda_{k-1}[c,\lambda_{k-1}[a,b]] = 0 $$ where
\begin{itemize}
  \item $\alpha = (-1)^{(|a|+k-1)(|c|+k-1)}$,
  \item $\beta = (-1)^{(|b|+k-1)(|a|+k-1)}$, and
  \item $\gamma = (-1)^{(|c|+k-1)(|b|+k-1)}$.
\end{itemize}

    \item (Antisymmetry)
$$\lambda_{k-1}[a,b] = (-1)^{|a||b|+1+(k-1)(|a|+|b|+1)} \lambda_{k-1}[b,a].$$

    \item (Product formula) $$\lambda_{k-1}[a \cdot b,c] =
    a \cdot \lambda_{k-1}[b,c]
    +  (-1)^{|a|\cdot|b|}b\cdot \lambda_{k-1}[a,c].$$
    \item (Commutation with homology suspension $\sigma_*$)
$$\sigma_*(\lambda_{k-1}[a,b]) =
    \lambda_{k-2}[\sigma_*(a),\sigma_*(b)].$$
    \item (Degree of the operation) The degree of $\lambda_{k-1}[a,b]$ is
$k-1 + |a| + |b|$.
\end{enumerate} In addition, it was proven that this pairing is
compatible with the Whitehead
product structure for the classical Hurewicz homomorphism via
the following commutative diagram

\[
\begin{CD}
\pi_{m+k}(X) \otimes \pi_{n+k}(X) @>{W_0 }>> \pi_{m+n + 2k-1}(X)  \\
   @VV{s_* \otimes s_*}V           @VV{s_*}V           \\
\pi_m(\Omega^kX) \otimes \pi_n(\Omega^kX) @>{ W_k }>> \pi_{m+n +
k-1}(\Omega^kX)  \\
   @VV{\phi \otimes \phi}V           @VV{\phi}V           \\
H_m(\Omega^kX) \otimes H_n(\Omega^kX) @>{\lambda_{k-1}}>> H_{m+n +
k-1}(\Omega^kX)
\end{CD}
\] for which

\begin{enumerate}
    \item the map $s_*$ is the natural isomorphism,
    \item the map $\phi$ is the classical Hurewicz homomorphism, and
    \item the map $W_k$ is the adjoint of the classical Whitehead
    product $$W_0: \pi_{m+k}(X) \otimes \pi_{n+k}X \to\ \pi_{m+n +
    2k-1}(X).$$
\end{enumerate}

\begin{thm} \label{thm: iterated loop poisson}
Assume that $k$ is greater than $1$.

\begin{enumerate}
\item If $k > 1$, the homology of
$\Omega^k({\rm Conf}^{G}(\mathbb H^2 \times \mathbb C^q,n))$, with any
field coefficients, is a graded Poisson algebra with Poisson bracket
given by the Browder operation $\lambda_{k-1}[-,-]$
for the homology of a $k$-fold loop space.

\smallskip

\item If $1< k <2q+1$, the homology of
$\Omega^k {\rm Conf}^{G}(\mathbb H^2 \times \mathbb C^q,n)$
with coefficients in a field $\bF$ of characteristic zero, is
the free Poisson algebra generated by elements
$$B_{i,j}^{\sigma}$$ of degree $2q+1 -k$ for $1 \leq j < i \leq n$, and
$\sigma$ in
$G$ modulo the ``infinitesimal Poisson surface
braid relations" given as follows:

      \begin{enumerate}
            \item If
            $\{i,j\} \cap \{s,t\} = \phi$, then
            $\lambda_{k-1}[B_{i,j}^{\sigma}, B_{s,t}^{\tau}] = 0$.
            \item If $ 1 \leq j < s < i \leq n$, then
            $\lambda_{k-1}[B^{\tau}_{i,j},B^{\tau \sigma^{-1} }_{i,s}+
B^{\sigma}_{s,j}] = 0$.
            \item If $ 1 \leq j < s < i \leq n $, then
            $\lambda_{k-1}[B^{\sigma}_{s,j},B^{\tau }_{i,j} + B^{\tau
\sigma^{-1}}_{i,s}] =0$.
            \item The antisymmetry relation, Jacobi identity, and product
            formula for a graded Poisson algebra are satisfied.
      \end{enumerate}

\smallskip

      \item There is a map $$E^2: {\rm Conf}^{G}(\mathbb H^2 \times \mathbb
C^q,n)
\to\ \Omega^2{\rm Conf}^{G}(\mathbb H^2 \times \mathbb C^{q+1},n) $$
which induces a homology isomorphism in degree $2q +1$. The associated loop map
$$\Omega(E^2):  \Omega {\rm Conf}^{G}(\mathbb H^2 \times \mathbb C^q,n)
\to\ \Omega^3{\rm Conf}^{G}(\mathbb H^2 \times \mathbb C^{q+1},n)$$ induces
an isomorphism on $H_{2q}(-;\Z)$. Furthermore, the image of the map
$\Omega(E^2)$
in homology is the subalgebra generated by the classes of degree $2q$.

\end{enumerate}
\end{thm}

\begin{proof}
Part (1) is a special case of results in \cite{CLM}, pages $215-217$.
Furthermore, if $ q > 1$, then the characteristic zero homology of
$\Omega^q(X)$ for a $q$-connected space $X$
is isomorphic to the graded symmetric algebra
generated by the image of the rational Hurewicz
homomorphism (as is well-known from work of Milnor-Moore ).

Part (2) now follows from the case $ q = 1 $ together
with the proof of Theorem \ref{thm:Lie algebras}. Namely,
the Browder operation $\lambda_{q-1}$ is precisely
the commutator for the underlying associative algebra given by the
homology of a $1$-fold loop space
$\Omega {\rm Conf}^{G}(\mathbb H^2 \times \mathbb C^q,n)$.
Furthermore, the homology suspension $\sigma_*$ satisfies
$ \sigma_*(\lambda_{k-2}[x,y]) = \lambda_{k-1}[\sigma_*(x),\sigma_*(y)]$
\cite[pages 215--217]{CLM}. Since the homology suspension when
restricted to the module of primitives $$\sigma_*:
\operatorname{Prim}H_n(\Omega^k(X);\mathbb Q)  \to\
\operatorname{Prim}H_{n+1}(\Omega^{k-1}(X);\mathbb Q)$$ is
an isomorphism when restricted to $k$-connected spaces $X$ for $k>2$,
it follows that the Poisson bracket relations follow at once from
the case $k=1$ as given in Theorem \ref{thm:Lie algebras}.

Part (3) requires a construction given as follows. Notice that
there is a map $$\Theta: {\rm Conf}^{G}(\mathbb H^2 \times \mathbb C^q,n)
\to\ \Pi_{I} S^{2q+1} $$
where the index set $I$ is given by $(i,j, \sigma, \tau)$ for $1 \leq j < i
\leq
n$ with $\sigma$, and $\tau$ elements of $G$. This map
is gotten by sending $(z_1, \dots, z_n)$ to $(\sigma(z_i) - \tau(z_j))/\rho$
for $1 \leq j < i \leq n$, with
$\sigma$, and $\tau$ in $G$, and $\rho = ||\sigma(z_i) - \tau(z_j)||$.

This map is a split monomorphism in homology. Thus after
suspending once, ${\rm Conf}^{G}(\mathbb H^2 \times \mathbb C^q,n)$ is homotopy
equivalent to a bouquet of spheres for any $ q \geq 0$.

Hence, there is a map
$\gamma:\Sigma^2 {\rm Conf}^{G}(\mathbb H^2 \times \mathbb C^q,n) \to\
\vee_I S^{2k+3}$ which
induces an isomorphism in homology in degree $2k+3$. There is an induced map
$$\sigma^2: \Sigma^2{\rm Conf}^{G}(\mathbb H^2 \times \mathbb C^q,n)  \to
{\rm Conf}^{G}(\mathbb H^2 \times \mathbb C^{q+1},n)$$ given by the composite
\[
\begin{CD}
\Sigma^2 {\rm Conf}^{G}(\mathbb H^2 \times \mathbb C^q,n)
@>{\gamma}>> \vee_I S^{2k+3} @>{\phi}>>
{\rm Conf}^{G}(\mathbb H^2 \times \mathbb C^{q+1},n)
\end{CD}
\] where $\phi: \vee_I S^{2k+3} \to\ {\rm Conf}^{G}(\mathbb H^2 \times
\mathbb C^{q+1},n) $ is the
natural map obtained from the Hurewicz homomorphism which induces
an isomorphism on the first non-vanishing homology group
of ${\rm Conf}^{G}(\mathbb H^2 \times \mathbb C^{q+1},n) $.

The adjoint of this composite
$$E^2: {\rm Conf}^{G}(\mathbb H^2 \times \mathbb C^q,n)  \to\ \Omega^2
{\rm Conf}^{G}(\mathbb H^2 \times \mathbb C^{q+1},n)$$ induces an
isomorphism on the first non-trivial homology group.  This last map
may be regarded as an analogue of the classical Freudenthal double
suspension map for which the spaces ${\rm Conf}^{G}(\mathbb H^2 \times
\mathbb C^q,n) $
are replaced by single odd dimensional spheres. Looping $E^2$ is given by
$\Omega(E^2): \Omega {\rm Conf}^{G}(\mathbb H^2 \times \mathbb C^q,n)
\to\ \Omega^3 {\rm Conf}^{G}(\mathbb H^2 \times \mathbb C^{q+1},n) $.

The theorem follows.
\end{proof}

Related structures for certain complements of
hyperplane arrangements are given in \cite{ccx}.

\end{document}